\documentclass[10pt]{article}

\usepackage{dirtytalk}

\usepackage{amsfonts}
\usepackage{amssymb}
\usepackage{amsmath}
\usepackage{enumerate}
\usepackage{amsthm}

\newcommand{\cf}{\mathrm{cf}}
\newcommand{\cof}{\mathrm{cof}}

\newcommand{\cov}{\mathrm{cov}}

\newcommand{\pp}{\mathrm{pp}}
\newcommand{\PP}{\mathrm{PP}}
\newcommand{\tcf}{\mathrm{tcf}}

\newtheorem{Th}{\bf THEOREM}[section]

\newtheorem{Lem}[Th]{\bf LEMMA}

\newtheorem{Pro}[Th]{\bf PROPOSITION}

\newtheorem{fact}[Th]{\bf FACT}

\newtheorem{Cor}[Th]{\bf COROLLARY}

\newtheorem{Obs}[Th]{\bf OBSERVATION}

\theoremstyle{definition} 
 
\theoremstyle{remark}

\theoremstyle{question}

\usepackage[english]{babel}
\makeindex

\title{$\mu$-CLUBS OF $P_\kappa (\lambda)$ : PARADISE ON EARTH}
\author{Pierre MATET}

\date{}

\begin{document}

\maketitle

\renewcommand{\thefootnote}{\arabic{footnote}} 	

\renewcommand{\thefootnote}{}                                



\vskip 0,7cm

\begin{abstract}  If $V = L$, and $\mu$, $\kappa$ and $\lambda$ are three infinite cardinals with $\mu = \cf (\mu) < \kappa = \cf (\kappa) \leq \lambda$, then, as shown in \cite{Heaven}, the $\mu$-club filters on $P_\kappa (\lambda)$ and $P_\kappa (\lambda^{< \kappa})$ are isomorphic if and only if $\cf (\lambda) \not= \mu$. Now in $L$, $\lambda^{< \kappa}$ equals $u (\kappa, \lambda)$ (the least size of a cofinal subset in $(P_\kappa (\lambda), \subseteq)$) equals $\lambda$ if $\cf (\lambda) \geq \kappa$, and $\lambda^+$ otherwise. We show that, in ZFC,  there are many triples $(\mu, \kappa, \lambda)$ for which ($u (\kappa, \lambda) > \lambda$ and) the $\mu$-club filters on $P_\kappa (\lambda)$ and $P_\kappa (u (\kappa, \lambda))$ are isomorphic.

\end{abstract}

\bigskip

\section{Introduction}

\bigskip


\subsection{Ask not what pcf theory can do for you ...}

\bigskip

True believers are a bit puzzled by the reception of Saharon Shelah's pcf theory. It has largely been seen as a tool, with most applications based on a kind of objects specific to the theory called \say{scales}. With each singular cardinal $\theta$ are associated scales, each with a length a regular cardinal $\pi$ with $\theta^+ \leq \pi \leq 2^\theta$. The scales considered are often those of length $\theta^+$, except for isolated, prizehunting attempts to find the value of $\aleph_\omega^{\aleph_0}$. Now the theory affirms that scales of length $\theta^+$ always exist, so the focus was on their properties, which led to a fine classification of their \say{points} (good, more-than-good, remarkably good, better, very good, etc.). This narrow focus may explain why interest started to wane after a while, and research on the subject slowed down. Let us observe that there is no reason why scales of different lengths should be all similar, and in fact scales of length $\theta^{++}$ (or of any length the successor of a regular cardinal) have nicer properties than their cousins of length $\theta^+$, a fact that we will take advantage of.

\bigskip

\subsection{Eine Weltanschauung}

\bigskip

For the faithful, pcf theory is not a tool, but a world to live in. They point out that the contents of The Book \cite{SheCA} do not match its title : the book can hardly be described as being about addition, multiplication and exponentiation of cardinal numbers. According to the devotees, the underlying message is clear : study pcf theory, not old style cardinal arithmetic. Pcf theory will tell you what you can prove, and will help you to prove it. You will thus obtain optimal results : results that are true in $L$, or more generally in $L$-like models (this is \say{paradise in heaven}), but also, for some values of the parameters involved, in ZFC (this is \say{paradise on earth}).

\bigskip

\subsection{Two-cardinal versions of $NS_\kappa \vert E^\kappa_\mu$}

\bigskip

Let us now describe our results in broad terms (for missing definitions see the next section). For two regular cardinals $\mu < \kappa$, $NS_\kappa \vert E^\kappa_\mu$ denotes the restriction of the nonstationary ideal on $\kappa$ to the set $E^\kappa_\mu$ of all limit ordinals $\delta$ less than $\kappa$ of cofinality $\mu$. It is identical with $N\mu$-$S_\kappa$, the ideal dual to the $\mu$-club filter on $\kappa$. Given a cardinal $\lambda \geq \kappa$, the respective generalizations of these two ideals to $P_\kappa (\lambda)$ (the restriction $NS_{\kappa, \lambda} \vert E^{\kappa, \lambda}_\mu$ of the nonstationary ideal on $P_\kappa (\lambda)$ to the set $E^{\kappa, \lambda}_\mu$ of all $a$ in $P_\kappa (\lambda)$ of uniform cofinality $\mu$, and the ideal $N\mu$-$S_{\kappa, \lambda}$ dual to the $\mu$-club filter on $P_\kappa (\lambda)$) do not necessarily coincide. Suppose for instance that $\kappa$ is weakly inaccessible, and $\lambda$ a strong limit cardinal of cofinality less than $\kappa$. Then $NS_{\kappa, \lambda}$ (and a fortiori $NS_{\kappa, \lambda} \vert E^{\kappa, \lambda}_\mu$) is a small ideal with few generators, since it is \cite{She02} a mere restriction of the noncofinal ideal $I_{\kappa, \lambda}$ on $P_\kappa (\lambda)$ (that has $u (\kappa, \lambda)$ many generators, where $u (\kappa, \lambda)$ denotes the least size of a cofinal subset in $(P_\kappa (\lambda), \subseteq)$), whereas $N\mu$-$S_{\kappa, \lambda}$ is a bigger ideal, since it is isomorphic to $N\mu$-$S_{\kappa, u (\kappa, \lambda)}$ and hence $cof (N\mu$-$S_{\kappa, \lambda}) = cof (N\mu$-$S_{\kappa, u (\kappa, \lambda)}) > u (\kappa, \lambda)$ \cite{MPS1}. On the other hand, if $\lambda < \kappa^{+ \omega}$, then $NS_{\kappa, \lambda} \vert E^{\kappa, \lambda}_\mu = N\mu$-$S_{\kappa, \lambda}$ \cite{Secret}.

\bigskip

\subsection{Pseudo-Kurepa families}

\bigskip

Assuming that $\cf (\lambda) < \kappa$, we will be looking for a $(\mu, \kappa, \lambda, u (\kappa, \lambda))$-shuttle, which is a function $\chi : P_\kappa (\lambda) \rightarrow P_\kappa (u (\kappa, \lambda))$ witnessing that $N\mu$-$S_{\kappa, \lambda}$ and $N\mu$-$S_{\kappa, u (\kappa, \lambda)}$ are isomorphic. A pseudo-Kurepa family for $P_\kappa (\lambda)$ is a subset $F$ of $P_\kappa (\lambda)$ with the property that $\vert F \cap P (b) \vert < \kappa$ for all $b \in P_\kappa (\lambda)$. The existence of a $(\mu, \kappa, \lambda, u (\kappa, \lambda))$-shuttle entails that of a pseudo-Kurepa family $F$ for $P_\kappa (\lambda)$ of size $u (\kappa, \lambda)$.
In the case when 
$\mu > \cf (\lambda)$, there is a partial converse : if there is a pseudo-Kurepa family $F$ for $P_\kappa (\lambda)$ with $F \subseteq P_{(\cf (\lambda))^+} (\lambda)$ and $\vert F \vert = u (\kappa, \lambda)$, then there exists a $(\mu, \kappa, \lambda, u (\kappa, \lambda))$-shuttle.

\bigskip

\subsection{The $\cov$ vs. $\pp$ problem}

\bigskip

Given four cardinals $\rho_1, \rho_2, \rho_3, \rho_4$ with $\rho_1 \geq \rho_2 \geq \rho_3 \geq \omega$ and $\rho_3 \geq \rho_4 \geq 2$, $\cov (\rho_1, \rho_2, \rho_3, \rho_4)$ denotes the least cardinality of any $X \subseteq  P_{\rho_2}(\rho_1)$ such that for any $a \in  P_{\rho_3}(\rho_1)$, there is $Q \in  P_{\rho_4}(X)$ with $a \subseteq \bigcup Q$. Thus $u (\kappa, \lambda) = \cov (\lambda, \kappa, \kappa, 2)$. The $\cov$ vs. $\pp$ problem asks whether everything under the sky (of pcf theory), and in particular covering numbers, can be expressed in terms of the $\pp$ functions. The specific version of the problem that we consider in this paper is whether $u (\kappa, \lambda) = \pp (\theta)$ for some $\theta$ with $\cf (\theta) < \kappa < \theta \leq \lambda$. In the case when $\pp (\theta)$ is regular, this means representing $u (\kappa, \lambda)$ as the length of a scale $\vec f$ on $\theta$. But we need this scale to be nice in the sense that all (or equivalently, almost all in the sense of the club filter on $u (\kappa, \lambda)$) points of $\vec f$ of cofinality $\kappa$ are good. From such a scale we can construct a pseudo-Kurepa family $F$ for $P_\kappa (\lambda)$ with the desired properties ($F \subseteq P_{(\cf (\theta))^+} (\lambda)$ and $\vert F \vert = u (\kappa, \lambda)$). But there is a glitch.  Whether points of $\vec f$ of cofinality $\kappa$ are good will often depend on both $\kappa$ and $\cf (\theta)$. Now this is a problem if, as in Fact 3.4 below, all you know is that there is a $\theta$ that works, with no clue given concerning the value of $\cf (\theta)$. So we are interested in situations when the $\theta$ that works can be computed, as in Fact 3.5.

\bigskip

\subsection{The long shadow of the RGCH Theorem}

\bigskip
 
A representative result of the present paper asserts that if $\lambda$ is a singular strong limit cardinal, and $\kappa$ a weakly inaccessible cardinal less than $\lambda$, then for each sufficiently large regular cardinal $\mu < \kappa$, $N\mu$-$S_{\kappa, \lambda}$ and $N\mu$-$S_{\kappa, u (\kappa, \lambda)}$ are isomorphic. A template for such statements is Shelah's Revised GCH Theorem \cite{SheRGCH} that affirms that given a cardinal $\lambda$, and a strong limit cardinal $\chi < \lambda$, $\cov (\lambda, \chi, \chi, \sigma) \leq \lambda$ for any sufficiently large cardinal $\sigma < \chi$. Since the threshold is not given, such results may seem of limited applicability.

\bigskip

\subsection{Sch\"{a}tzchen}

\bigskip

They are however of great interest from a theoretical point of view, since they show that ZFC is able to decide by itself some serious issues, and this is in a nice way. In this sense, ZFC is more complete than previously feared. This is why some devotees are of the opinion that results of this type should be counted among the \say{treasures} that, according to Shelah \cite{SheCantor}, \say{are waiting for you}.

\bigskip

\section{Shuttles}

\bigskip

\subsection{Ideals}

\bigskip

For a set $A$ and a cardinal $\chi$, we set $P_\chi (A) = \{ a\subseteq A :  \vert a \vert < \chi\}$ and $ [A]^\chi = \{ a\subseteq A :  \vert a \vert  = \chi\}$.

\medskip

Let $X$ be an infinite set. An {\it ideal  on} $X$ is a nonempty collection $J$ of subsets of $X$ such that 
 \begin{itemize}
\item $X \notin J$.
\item $P(A)\subseteq J$ for all $A\in J$.  
\item $A\cup B\in J$ whenever $A, B \in J$. 
\end{itemize}

\medskip

Given an ideal $J$ on $X$, we let $J^+ = P(X) \setminus J$, $J^\ast = \{ A\subseteq X : X \setminus A\in J\}$, and $J \vert A = \{ B\subseteq X : B \cap A\in J\}$  for each $A\in J^+$.  
For a cardinal   $\rho,J$ is {\it $\rho$-complete} if  $\bigcup Q \in J$ for every $Q \subseteq J$ with $\vert Q \vert< \rho$. ${\rm cof} (J)$ denotes the least cardinality of any $Q \subseteq J$ such that $J = \bigcup_{A\in Q} P(A)$. 

\medskip

Let $X$ and $Y$ be two infinite sets, and $J$ be an ideal  on $X$. Given $f : X \rightarrow Y$,  we let $f (J) = \{ B \subseteq Y : f^{-1} (B) \in J\}$.   

\medskip

\begin{fact} {\rm (\cite{Secret})} Suppose that $f^{-1} (\emptyset) \in J$. Then 
\begin{enumerate}[\rm (i)]
\item $f (J)$ is an ideal on $Y$.
\item If $J$ is $\rho$-complete, then so is $f (J)$.
\item $(f (J))^\ast = \{C \subseteq Y : f^{-1} (C) \in J^\ast\}$.
\item $\{ f`` D : D \in J^\ast\} \subseteq (f (J))^\ast$.
\item $\cof (f (J)) \leq \cof (J)$.
\end{enumerate}
\end{fact} 

\medskip

For $i = 1, 2$, let $X_i$  be an infinite set, and $K_i$ be an ideal  on $X_i$. We say that $K_1$ is {\it isomorphic to} $K_2$ if there are $W_1 \in K_1^\ast$, $W_2 \in K_2^\ast$ and a bijection $k : W_1 \rightarrow W_2$ such that 

\centerline{$K_1^\ast \cap P (W_1) = \{D \subseteq W_1 : k``D \in K_2^\ast\}$.}
\medskip




\begin{Obs} 
\begin{enumerate}[\rm (i)]
\item Suppose that $K_1$ is isomorphic to $K_2$. Then $K_2$ is isomorphic to $K_1$.
\item Suppose that $K_1$ is isomorphic to $K_2$, and $K_2$ isomorphic to $K_3$, where $K_3$ is an ideal on some infinite set $X_3$. Then $K_1$ is isomorphic to $K_3$.
\item Suppose that $K_1$ is isomorphic to $K_2$, and let $Y_1 \in K_1^+$. Then $K_1 \vert Y_1$ is isomorphic to $K_2 \vert T_2$ for some $T_2 \in K_2^+$.
\item Suppose that $W_1 \in K_1^\ast$, $W_2 \in K_2^\ast$ and $k : W_1 \rightarrow W_2$ are such that $K_1^\ast \cap P (W_1) = \{D \subseteq W_1 : k``D \in K_2^\ast\}$. Then $k (K_1) = K_2$.

\end{enumerate}
\end{Obs} 

{\bf Proof.} (i) : For a proof see e.g. \cite{Menas}.

\medskip

(ii) :  We may find $W_1 \in K_1^\ast$, $W_2, Y_2 \in K_2^\ast$, $Y_3 \in K_3^\ast$ and bijections $k : W_1 \rightarrow W_2$ and $h : Y_1 \rightarrow Y_2$ such that $K_1^\ast \cap P (W_1) = \{D \subseteq W_1 : k``D \in K_2^\ast\}$ and $K_2^\ast \cap P (Y_2) = \{C \subseteq Y_2 : h``C \in K_3^\ast\}$. Set $T_2 = W_2 \cap Y_2$, $T_1 = k^{- 1} (T_2)$, $T_3 = h``T_2$, and define $g : T_1 \rightarrow T_3$ by $g (x) = h (k (x))$. Then, as is readily checked, $T_i \in K_i^\ast$ for $i = 1, 2, 3$, $g$ is a bijection, and $K_1^\ast \cap T_1 = \{B \subseteq T_1 : g``B \in K_3^\ast\}$.

\medskip

(iii) : We may find $W_1 \in K_1^\ast$, $W_2 \in K_2^\ast$ and a bijection $k : W_1 \rightarrow W_2$ such that $K_1^\ast \cap P (W_1) = \{D \subseteq W_1 : k``D \in K_2^\ast\}$. Set $T_1 = Y_1 \cap W_1$, $T_2 = k`T_1$, and define $g : T_1 \rightarrow T_2$ by $g (x) = k (x)$. Then clearly $T_1 \in (K_1 \vert Y_1)^\ast$, and moreover $g$ is a bijection.
\medskip

{\bf Claim 1.}  $T_2 \in K_2^+$.      

\medskip

{\bf Proof of Claim 1.} Suppose otherwise. Set $B_2 = W_2 \setminus T_2$. Then $B_2 \in K_2^\ast$, so $k^{- 1} (B_2) \in K_1^\ast$, and therefore $T_1 \cap k^{- 1} (B_2) \in K_1^+$. But it is easy to see that $T_1 \cap k^{- 1} (B_2) = \emptyset$. This contradiction completes the proof of the claim. 

\medskip 

It remains to show that $(K_1 \vert Y_1)^\ast  \cap P (T_1) = \{C \subseteq T_1 : g``C \in (K_2 \vert T_2)^\ast\}$.

\medskip

{\bf Claim 2.}  Let $C \in (K_1 \vert Y_1)^\ast \cap P (T_1)$.  Then $g``C \in (K_2 \vert T_2)^\ast$.    

\medskip

{\bf Proof of Claim 2.} Select $H \in K_1^\ast$ with $H \cap T_1 \subseteq C$. Then $k``(H \cap W_1) \in K_2^\ast$. Now

\centerline{$k``(H \cap W_1) \cap T_2 = k``(H \cap W_1) \cap k`T_1 = k``(H \cap T_1) = g``(H \cap T_1) \subseteq g``C$,}
 
which completes the proof of the claim. 

\medskip 

{\bf Claim 3.}  Let $C \subseteq T_1$ be such that $g``C \in (K_2 \vert T_2)^\ast$. Then $C \in (K_1 \vert Y_1)^\ast$.    

\medskip

{\bf Proof of Claim 3.} Pick $G \in K_2^\ast$ so that $G \cap T_2 \subseteq g``C$. Then 

\centerline{$k^{- 1} (G \cap W_2) \cap T_1 = k^{- 1} (G \cap W_2) \cap k^{- 1} (T_2) = k^{- 1} (G \cap T_2) = g^{- 1} (G \cap T_2) \subseteq C$.}
 
Since $k^{- 1} (G \cap W_2) \in K_1^\ast$ and $K_1 \vert T_1 = K_1 \vert Y_1$, it follows that $C \in (K_1 \vert Y_1)^\ast$, which completes the proof of the claim and that of (iii).  

\medskip

(iv) : Given $T \subseteq X_2^\ast$, $T \in K_2^\ast$ if and only if $T \cap ran (k) \in K_2^\ast$ if and only if $k^{- 1} (T \cap ran (k)) \in K_1^\ast$ if and only if $k^{- 1} (T) \cap k^{- 1} (ran (k)) \in K_1^\ast$ if and only if $k^{- 1} (T) \cap W_1 \in K_1^\ast$ if and only if $k^{- 1} (T) \in K_1^\ast$ if and only if $T \in (k (K_1))^\ast$.
\hfill$\square$  

\begin{fact}  {\rm(\cite{Heaven})} The following are equivalent :
\begin{enumerate}[\rm (i)]
\item  $f (K_1) = K_2$ for some one-to-one $f : X_1 \rightarrow X_2$.
\item  There are $W_2 \in K_2^\ast$ and a bijection $k : X_1 \rightarrow W_2$ such that $K_1^\ast = \{D \subseteq X_1 : k``D \in K_2^\ast\}$.
\end{enumerate}
\end{fact}

\medskip

Thus if $f (K_1) = K_2$ for some one-to-one $f : X_1 \rightarrow X_2$, then $K_1$ and $K_2$ are isomorphic.

\bigskip

\subsection{$\mu$-clubs}

\bigskip

Given a regular uncountable cardinal $\kappa$ and a cardinal $\sigma \geq \chi$, $I_{\kappa, \sigma}$ and $NS_{\kappa, \sigma}$ denote, respectively,  the noncofinal ideal on $P_\kappa (\sigma)$ and the nonstationary ideal on $P_\kappa (\sigma)$.  
 
 \medskip
 
 An ideal $J$ on $P_\kappa (\sigma)$ is {\it fine} if $I_{\kappa,\sigma} \subseteq J$. 

\medskip

For a regular cardinal $\mu < \kappa$, a subset $A$ of $P_\kappa (\sigma)$ is $\mu${\it -closed} if $\bigcup_{i < \mu} a_i \in A$ for every increasing sequence $\langle a_i : i < \mu \rangle$ in $(A, \subset)$. A subset $C$ of $P_\kappa (\sigma)$ is a $\mu${\it -club} if it is a $\mu$-closed, cofinal subset of $P_\kappa (\sigma)$. 

\medskip

We let  $N\mu$-$S_{\kappa,\sigma}$ be the set of all $B \subseteq P_\kappa (\sigma)$ such that $B \cap C = \emptyset$ for some $\mu$-club $C\subseteq P_\kappa (\sigma)$. 

\medskip

\begin{fact}
{\rm(\cite{Conc})}  $N\mu$-$S_{\kappa,\sigma}$ is a normal ideal on $P_\kappa (\sigma)$.
\end{fact} 


\medskip

Let $\pi \geq \sigma$ be a cardinal. A $(\mu, \kappa, \sigma, \pi)$-{\it shuttle} is a function $\chi : P_\kappa (\sigma) \rightarrow P_\kappa (\pi)$ such that
\begin{enumerate}[\rm (a)]
\item $\chi (a) \cap \sigma = a$ for all $a \in P_\kappa (\sigma)$.
\item $\chi (a) \subseteq \chi (b)$ whenever $a, b \in P_\kappa (\sigma)$ are such that $a \subseteq b$.
\item $ran (\chi) \in I^+_{\kappa, \pi}$.
\item $\chi (\bigcup_{i < \mu} a_i) \subseteq \bigcup_{i < \mu} \chi (a_i)$ for any increasing sequence $\langle a_i : i < \mu \rangle$ in $(P_\kappa (\sigma), \subset)$.
\end{enumerate}

\medskip

Notice that in case $\pi = \sigma$, the identity function on $P_\kappa (\sigma)$ is a $(\mu, \kappa, \sigma, \pi)$-shuttle, and it is the only one.

\medskip

\begin{fact} {\rm(\cite{Heaven})} Let $\chi : P_\kappa (\sigma) \rightarrow P_\kappa (\pi)$ be a $(\mu, \kappa, \sigma, \pi)$-shuttle. Then the following hold :
\begin{enumerate}[\rm (i)]
\item  $\chi$ is one-to-one.
\item $a \subseteq b$ whenever $a, b \in P_\kappa (\sigma)$ are such that $\chi (a) \subseteq \chi (b)$.
\item $\chi (a) \subset \chi (b)$ whenever $a, b \in P_\kappa (\sigma)$ are such that $a \subset b$.
\item $\chi`` A \in I^+_{\kappa, \pi}$ for all $A \in I^+_{\kappa, \sigma}$.
\item $\chi (\bigcup_{i < \mu} a_i) = \bigcup_{i < \mu} \chi (a_i)$ for any increasing sequence $\langle a_i : i < \mu \rangle$ in $(P_\kappa (\sigma), \subset)$.
\item $\chi (I_{\kappa, \sigma}) = I_{\kappa, \pi} \vert ran (\chi)$.
\end{enumerate}
\end{fact}

\medskip

It follows from Facts 2.3 and 2.5 that the existence of a $(\mu, \kappa, \sigma, \pi)$-shuttle entails that $I_{\kappa, \sigma}$ is isomorphic to a restriction of $I_{\kappa, \pi}$.

\medskip

Given two uncountable cardinals $\tau = \cf (\tau) \leq \rho$, we let $u (\tau, \rho) =$ the least size of a cofinal set in $(P_\tau (\rho), \subseteq)$.

\medskip

\begin{fact} {\rm(\cite{Heaven})}  Suppose that there exists a $(\mu, \kappa, \sigma, \pi)$-shuttle $\chi$. Then 
\begin{enumerate}[\rm (i)]
\item $\pi \leq u (\kappa, \sigma) = u (\kappa, \pi)$.
\item $\chi (N\mu$-$S_{\kappa,\sigma}) = N\mu$-$S_{\kappa,\pi}$. 
\item $N\mu$-$S_{\kappa,\sigma}$ and $N\mu$-$S_{\kappa,\pi}$ are isomorphic.
\item $p_\pi (N\mu$-$S_{\kappa,\pi}) = N\mu$-$S_{\kappa,\lambda}$, where $p_\pi : P_\kappa (\pi) \rightarrow P_\kappa (\lambda)$ is the projection defined by $p_\pi (x) = x \cap \lambda$.
\end{enumerate}
\end{fact}

\bigskip

\subsection{${\mathcal A}_{\kappa,\lambda} (\tau,\pi)$}

\bigskip

Let $\tau, \kappa, \lambda, \pi$ be four cardinals with $2 \leq \tau$ and $\max \{ \tau, \omega \} \leq \kappa = \cf (\kappa) \leq \lambda \leq \tau$. A $(\tau,\lambda,\pi)${\it -sequence} is a one-to-one sequence $\vec y = \langle y_\beta : \beta < \pi \rangle$ of elements of $P_\tau (\lambda)$  with the property that $y_\beta = \{\beta\}$  for every $\beta < \lambda$.

Given a $(\tau,\lambda,\pi)$-sequence $\vec y$, define $f_{\vec y } : P_\kappa (\lambda) \rightarrow P (\pi)$ by $f_{\vec y } (a) = \{\beta < \pi : y_\beta \subseteq a\}$, and let

\centerline{$A_\kappa (\vec y) = \{x \in P_\kappa (\pi) : f_{\vec y } (x \cap \lambda) \subseteq x\}$}

and 
 
\centerline{$\Delta_\kappa (\vec y) = \{x \in P_\kappa (\pi) : \forall \delta \in x (y_\delta \subseteq x)\}$.} 

\medskip

\begin{fact} {\rm (\cite{Secret}, \cite{Heaven})}
\begin{enumerate}[\rm (i)]
\item $a =  f_{\vec y} (a) \cap \lambda$.
\item $A_\kappa (\vec y) \cap \Delta_\kappa (\vec y) = \{x \in P_\kappa (\pi) : f_{\vec y} (x \cap \lambda) = x\} = ran ( f_{\vec y})$.
\item $\Delta_\kappa (\vec y) \in NS^\ast_{\kappa, \pi}$. 
\item  $A_\kappa (\vec y) \in I^+_{\kappa, \pi}$ if and only if $\vert f_{\vec y} (a) \vert < \kappa$ for every $a \in P_\kappa (\lambda)$
\end{enumerate}
\end{fact}

\medskip

An ${\mathcal A}_{\kappa,\lambda} (\tau,\pi)${\it -sequence} is a a $(\tau,\lambda,\pi)$-sequence $\vec y$ with the property that $\vert f_{\vec y} (a) \vert < \kappa$ for every $a \in P_\kappa (\lambda)$.

${\mathcal A}_{\kappa,\lambda} (\tau,\pi)$ asserts the existence of an ${\mathcal A}_{\kappa,\lambda} (\tau,\pi)$-sequence.

\medskip

\begin{fact} {\rm(\cite{Weaksat})}
\begin{enumerate}[\rm (i)]
\item $\mathcal{A}_{\kappa,\lambda} (2,\lambda)$ holds.
\item If $\mathcal{A}_{\kappa,\lambda} (\tau,\pi)$ holds, then $\pi \leq \cov (\lambda, \kappa, \tau, 2)$.
\end{enumerate}  
\end{fact}

\medskip

In particular, $\mathcal{A}_{\kappa,\lambda} (\kappa,\pi)$ entails that $\pi \leq u (\kappa, \lambda)$.

\medskip

\begin{Obs} Suppose that ${\mathcal A}_{\kappa,\lambda} (\tau,\pi)$ holds. Then for any cardinal $\rho \geq \lambda$, ${\mathcal A}_{\kappa,\rho} (\tau,\pi)$ holds.
\end{Obs}

\begin{fact} {\rm(\cite{Heaven})}  Let $\vec y = \langle y_\delta : \delta < \pi \rangle$ be an ${\mathcal A}_{\kappa,\lambda} (\kappa,\pi)$-sequence. Then the following hold : 
\begin{enumerate}[\rm (i)]
\item $ran ( f_{\vec y}) \in I_{\kappa, \pi}^+$.
\item $f_{\vec y} $ is an isomorphism from $(P_\kappa (\lambda), \subset)$ to $(ran ( f_{\vec y}), \subset)$.
\item $I_{\kappa, \pi} \vert ran ( f_{\vec y}) = f_{\vec y} (I_{\kappa, \lambda})$. 
\item $I_{\kappa, \lambda}$ and $I_{\kappa, \pi} \vert ran (f_{\vec y})$ are isomorphic.
\end{enumerate}
\end{fact}

\begin{Obs}  The following are equivalent :
\begin{enumerate}[\rm (i)]
\item ${\mathcal A}_{\kappa,\lambda} (\kappa, \pi)$ holds.
\item For some $X \in I_{\kappa, \pi}^+$, $I_{\kappa, \lambda}$ and $I_{\kappa, \pi} \vert X$ are isomorphic.
\end{enumerate}
\end{Obs}

{\bf Proof.}

\hskip0,2cm  (i) $\rightarrow$ (ii) : By Fact 2.10 (iv).

\medskip

\hskip0,2cm  (ii) $\rightarrow$ (i) : Suppose that $X \in I_{\kappa, \pi}^+$ is such that $I_{\kappa, \lambda}$ and $I_{\kappa, \pi} \vert X$ are isomorphic. There must be $W_1 \in (I_{\kappa, \pi} \vert X)^\ast$, $W_2 \in I_{\kappa, \lambda}^\ast$ and a bijection $k : W_1 \rightarrow W_2$ such that 

\centerline{$(I_{\kappa, \pi} \vert X)^\ast \cap P (W_1) = \{D \subseteq W_1 : k``D \in I_{\kappa, \lambda}^\ast\}$.}

For each $\gamma < \pi$, select $a_\gamma \in P_\kappa (\lambda)$ such that

\centerline{$\{b \in P_\kappa (\lambda) : a_\gamma \subseteq b \} \subseteq k``\{ x \in X \cap W_1 : \gamma \in x \}$.}

Now assume to the contrary that $\vert \bigcup_{\gamma \in e} a_\gamma \vert < \kappa$ for some $e \in [\pi]^\kappa$. For each $\delta \in e$, pick $x_\delta \in X \cap W_1$ so that $\delta \in x_\delta$ and $k (x_\delta) = \bigcup_{\gamma \in e} a_\gamma$. Since $k$ is one-to-one, we may find $x$ such that $x_\delta = x$ for all $\delta \in e$. But then $e \subseteq x$, which yields the desired contradiction. Thus $\langle a_\gamma : \gamma < \pi \rangle$ is an ${\mathcal A}_{\kappa,\lambda} (\kappa,\pi)$-sequence.
\hfill$\square$  

\bigskip

\section{Large $\mu$}

\bigskip

\subsection{Shuttles from pseudo-Kurepa families}

\bigskip
\begin{fact} {\rm(\cite{Heaven})} Let  $\vec y = \langle y_\beta : \beta < \pi \rangle$ be an ${\mathcal A}_{\kappa,\lambda} (\mu,\pi)$-sequence. Then $f_{\vec y}$ is a $(\mu, \kappa, \lambda, \pi)$-shuttle.
\end{fact}

\medskip

Notice that by Observation 2.11 and Facts 2.3 and 2.5, ${\mathcal A}_{\kappa,\lambda} (\kappa,\pi)$ follows from the existence of a $(\mu, \kappa, \lambda, \pi)$-shuttle. As shown by Todorcevic (see \cite{Weaksat}), there may exist (relative to large large cardinals) situations when $\kappa$ is a successor cardinal, $\cf (\lambda) < \kappa$ and ${\mathcal A}_{\kappa,\lambda} (\kappa,\lambda^+)$ fails.

\medskip

\begin{Obs} Let $\mu$, $\kappa$, $\lambda$, $\sigma$ be four infinite cardinals such that $\cf (\lambda) < \mu = \cf (\mu) < \kappa = \cf (\kappa) < \lambda \leq \sigma \leq u (\kappa, \lambda)$. Suppose that ${\mathcal A}_{\kappa,\lambda} (\mu, u (\kappa, \lambda))$ holds. Then $N\mu$-$S_{\kappa,\lambda}$ and $N\mu$-$S_{\kappa,\sigma}$ are isomorphic.
\end{Obs}

{\bf Proof.} 

{\bf Claim.}  Let $\rho$ be a cardinal with $\lambda \leq \rho \leq u (\kappa, \lambda)$. Then $N\mu$-$S_{\kappa,\rho}$ and $N\mu$-$S_{\kappa,u (\kappa, \lambda)}$ are isomorphic.
\medskip

{\bf Proof of the claim.} By Observation 2.9, ${\mathcal A}_{\kappa,\rho} (\mu, u (\kappa, \lambda))$ holds, so by Facts 2.6 and 3.1, $N\mu$-$S_{\kappa,\sigma}$ and $N\mu$-$S_{\kappa,u (\kappa, \lambda)}$ are isomorphic, which completes the proof of the claim. 

\medskip 
By the claim, $N\mu$-$S_{\kappa,\lambda}$ and $N\mu$-$S_{\kappa,u (\kappa, \lambda)}$ are isomorphic, and so are $N\mu$-$S_{\kappa,\sigma}$ and $N\mu$-$S_{\kappa,u (\kappa, \lambda)}$. Now use Observation 2.2 (ii).
\hfill$\square$

\bigskip

\subsection{Scales with good points}

\bigskip

One way to obtain ${\mathcal A}_{\kappa,\lambda} (\mu, u (\kappa, \lambda))$ is to use scales.

\medskip

Let $A$ be an infinite set of regular infinite cardinals such that $\vert A \vert  < \min A$, and $I$ be an ideal on $A$. For $f, g \in {}^A On$, let $f <_I g$ if $\{ a \in A : f(a) < g(a)\} \in I^\ast$,  and  $f\leq_I g$ if $\{ a\in A : f(a) \leq g(a)\} \in I^\ast$.

We let $\prod A =  \prod_{a\in A} a$.

Given a regular cardinal $\pi$, we say that $\tcf (\prod A/I) = \pi$ if there exists an increasing, cofinal sequence $\vec{f} = \langle f_\xi : \xi < \pi\rangle$ in $(\prod A, <_I)$. 
Such a sequence is said to be a {\it scale of length} $\pi$ for $\sup (A)$.

\medskip

For a singular cardinal $\theta$, we let 

\centerline{$\pp (\theta) = \sup \{ \pi : \pi$ is the length of some scale for $\theta \}$.}

\medskip

Suppose that $\tcf (\prod A/I) = \pi$, where $A$ is an infinite set of regular infinite cardinals such that $\vert A \vert  < \min (A)$, and $I$ an ideal on $A$, and that $\vec{f} = \langle f_\xi : \xi < \pi \rangle$ is an $<_I$-increasing, cofinal sequence in $(\prod A, <_I)$. An infinite limit ordinal $\delta < \pi$  is a {\it good point} for $\vec f$ if we may find a cofinal subset $X$ of $\delta$,  and $Z_\xi\in I$ for $\xi\in X$ such that $f_\beta (i) < f_\xi (i)$  whenever $\beta < \xi$  are in $X$ and $i \in \nu \setminus (Z_\beta \cup Z_\xi)$. 

\medskip

For an infinite cardinal $\tau$, we let $\rho (\tau)$ denote the largest limit cardinal less than or equal to $\tau$.

\medskip

\begin{fact}  {\rm(\cite{Good}, \cite{Weaksat})} Let $\theta, \pi, A$ and $I$ be such that
\begin{itemize}
\item $\theta$ and $\pi$ are two uncountable cardinals with $\theta < \pi$.
\item $A$ is a set of regular cardinals smaller than $\theta$ with supremum $\theta$.
\item $\vert A \vert < \min A$.
\item $I$ is an ideal on $A$ with $\{A \cap a : a \in A \} \subseteq I$.
\item  $\pi = \tcf (\prod A /I )$.
\end{itemize}
Further let  $\vec{f} = \langle f_\alpha : \alpha < \pi \rangle$ be an increasing, cofinal sequence in $(\prod A, <_I)$. Then the following hold :
\begin{enumerate}[\rm (i)]
\item There is a closed unbounded subset $C_{\vec{f}}$ of $\pi$, consisting of infinite limit ordinals, with the property that any $\delta$ in $C_{\vec{f}}$ satisfying one of the following conditions, where $\rho (\cf (\delta))$ is abbreviated as $\rho$, is a good point for $\vec f$ :
\begin{enumerate}[\rm (a)]
\item $(\max \{ \rho, \vert A \vert \})^{+3} < \cf (\delta)$.
\item $\rho^{\vert A \vert} < \cf (\delta)$.
\item $\vert A \vert < \cf (\rho)$.
\item $\vert A \vert < \rho$  and $I$ is $(\cf (\rho))^+$-complete.
\item $\cf (\rho) \leq \vert A \vert < \rho$ and $\pp (\rho) < \cf (\delta)$. 
\end{enumerate}
\item Let $\kappa$ be a successor cardinal with $\vert A \vert < \kappa < \theta$. Suppose that there is a closed unbounded subset $C$ of $\pi$ such that every $\delta \in C$ of cofinality $\kappa$ is a good point for ${\vec f}$. Then ${\mathcal A}_{\kappa,\theta} (\vert A \vert ^+, \pi)$ holds.
\item Let $\kappa$ be a weakly inaccessible cardinal with $\vert A \vert < \kappa < \theta$. Then ${\mathcal A}_{\kappa,\theta} (\vert A \vert ^+, \pi)$ holds.
\end{enumerate}
\end{fact}

\bigskip

\subsection{Non fixed points}

\bigskip

We now concentrate on those singular cardinals that are not fixed points of the aleph function, the reason being that they are much easier to handle. We first recall some definitions.

\medskip

Given three infinite cardinals $\sigma$, $\tau$ and $\theta$ with $\sigma \leq \cf (\theta) < \tau < \theta$, we let $\PP_{\Gamma (\tau, \sigma)} (\theta)$ be the collection of all cardinals $\pi$ such that $\pi = \tcf( \prod A /I )$ for some $A$ and $I$ such that
\begin{itemize}
\item $A$ is a set of regular cardinals smaller than $\theta$.
\item $\sup A = \theta$.
\item $\vert A \vert < \min \{\tau, \min A\}$.
\item $I$ is a $\sigma$-complete ideal on $A$ such that $\{A \cap a : a \in A \} \subseteq I$.
\end{itemize}

We let $\pp_{\Gamma (\tau, \sigma)} (\theta) = \sup \PP_{\Gamma (\tau, \sigma)} (\theta)$.

\medskip

Considering the special case when $\sigma = \omega$, we put $\PP_{< \tau} (\theta) = \PP_{\Gamma (\tau, \omega)} (\theta)$ and $\pp_{< \tau} (\theta) = \pp_{\Gamma (\tau, \omega)} (\theta)$.

\medskip

We let $\PP_\nu (\theta) = \PP_{< \nu^+} (\theta)$ and $\pp_\nu (\theta) = \sup \PP_\nu (\theta)$ for each cardinal $\nu$ with $\cf (\theta) \leq \nu < \theta$.

\medskip

Thus $\pp (\theta) = \pp_{\cf (\theta)} (\theta)$.

\medskip

For a cardinal $k$, $FP(k)$ denotes the least fixed point of the aleph function greater than $k$. 

\medskip

\begin{fact} {\rm (\cite{pcf})} Let $\kappa$ be an infinite successor cardinal, and $\lambda$ be a cardinal with $\kappa < \lambda < \min\{FP(\kappa), u(\kappa, \lambda)\}$. Then there is a (unique) cardinal $\theta (\kappa, \lambda)$ such that
 \begin{itemize}
 \item $\cf(\theta (\kappa, \lambda)) < \kappa < \theta (\kappa, \lambda) \leq \lambda$.
\item $u(\kappa, \lambda) = u (\kappa, \theta (\kappa, \lambda)) = \pp (\theta (\kappa, \lambda)) = \pp_{\Gamma((\cf(\theta (\kappa, \lambda))^+, \cf(\theta (\kappa, \lambda)))} (\theta (\kappa, \lambda)$.
\end{itemize}
\end{fact}

\begin{fact} {\rm (\cite{Fixed})} Let $\theta$ be a singular cardinal that is not a fixed point of the aleph function. Then the following hold : 
\begin{enumerate}[\rm (i)]
\item $\pp (\theta)$ is a successor cardinal.
\item $\pp (\theta) = \pp_{< \theta} (\theta) = \pp_{\Gamma((\cf(\theta)^+, \cf(\theta))} (\theta)$.
\item For any large enough regular cardinal $\kappa < \theta$, the following hold :
\begin{enumerate}[\rm (a)]
\item $u (\kappa, \theta) = \pp (\theta)$ and $\theta (\kappa, \theta) = \theta$.
\item $u((\cf (\theta))^+, \theta) = \max \{u((\cf (\theta))^+, \kappa), u(\kappa, \theta) \}$.
\item $u (\kappa, \rho) < \theta$ for every cardinal $\rho$ with $\kappa \leq \rho < \theta$.
\end{enumerate}
\end{enumerate}
\end{fact}

\begin{Th} Let $\theta$ be a singular cardinal that is not a fixed point of the aleph function. Then for any sufficiently large regular cardinal $\kappa < \theta$ such that either $(\rho (\kappa))^{+ 3} < \kappa$, or $\cf (\rho (\kappa)) \not= \cf (\theta)$, $u(\kappa, \theta) = \pp (\theta)$, and moreover ${\mathcal A}_{\kappa,\theta} ((\cf (\theta)) ^+, \pp (\theta))$ holds. 
\end{Th}

{\bf Proof.} By Fact 3.5, $\pp (\theta)$ is a successor cardinal, so there must be $A$ and $I$ such that
\begin{itemize}
\item $A$ is a set of regular cardinals smaller than $\theta$ with supremum $\theta$.
\item $\vert A \vert = \cf (\theta) < \min A$.
\item $I$ is a $\cf (\theta)$-complete ideal on $A$ with $\{A \cap a : a \in A \} \subseteq I$.
\item  $\pp (\theta) = \tcf (\prod A /I )$.
\end{itemize}
Select an increasing, cofinal sequence $\vec{f} = \langle f_\alpha : \alpha < \pp (\theta) \rangle$ in $(\prod A, <_I)$.

\medskip

By Fact 3.5, we may find $\eta$ with $\cf (\theta) < \eta < \theta$ such that for any regular cardinal $\kappa$ with $\eta \leq \kappa < \theta$, the following hold :
\begin{itemize}
\item $\kappa$ is not a fixed point of the aleph function (so it is a successor cardinal).
\item $u (\kappa, \theta) > \theta = \theta (\kappa, \theta)$.
\end{itemize}
Fix a regular cardinal $\kappa$ such that $\eta \leq \kappa < \theta$, $\cf (\theta) \leq \rho (\kappa)$ and either $(\rho (\kappa))^{+ 3} < \kappa$, or $\cf (\rho (\kappa)) \not= \cf (\theta)$. By Fact 3.3 ((i) and (ii)), there is a closed unbounded subset $C$ of $\pp (\theta)$ such that every $\delta \in C$ of cofinality $\kappa$ is a good point for ${\vec f}$, and hence ${\mathcal A}_{\kappa,\theta} (\vert A \vert ^+, \pp (\theta))$ holds.
\hfill$\square$

\begin{Th} Let $\kappa$ be an infinite successor cardinal, and $\lambda$ a cardinal greater than $\kappa$. Suppose that either $(\rho (\kappa))^{+ 3} < \kappa$ and $\lambda < \kappa^{+ (\rho (\kappa))^+}$, or $\lambda < \kappa^{+ \cf (\rho (\kappa))}$. Then ${\mathcal A}_{\kappa, \lambda} ((\cf (\theta (\kappa, \lambda)))^+, u (\kappa, \lambda))$ holds. 
\end{Th}

{\bf Proof.} Notice that if $\lambda < \kappa^{+ (\rho (\kappa))^+}$ (respectively, $\lambda < \kappa^{+ \cf (\rho (\kappa))}$), then $\cf (\theta (\kappa, \lambda)) \leq \rho (\kappa)$, (respectively, $\cf (\theta (\kappa, \lambda)) < \cf (\rho (\kappa))$), and proceed as in the proof of Theorem 3.6.
\hfill$\square$

\medskip

Thus in particular, ${\mathcal A}_{\kappa, \lambda} (\omega_1, u (\kappa, \lambda))$ holds whenever $\kappa = \omega_n$ with $4 \leq n < \omega$, and $\kappa \leq \lambda < \omega_{\omega_1}$.

\bigskip

\subsection{Non fixed points 2}

\bigskip

Theorems 3.6 and 3.7 deal with successor cardinals $\kappa$. This subsection is concerned with the case when $\kappa$ is weakly inaccessible.

\medskip

\begin{fact} {\rm (\cite{pcf})} Let $\kappa$ be a weakly inaccessible cardinal, and $\lambda$ be a cardinal with $\kappa^{+(\kappa \cdot n)} \leq \lambda < \kappa^{+(\kappa \cdot (n+1))}$, where $n < \omega$. Then the following hold :
\begin{enumerate}[\rm (i)]
\item Suppose that $\lambda = \kappa^{+(\kappa \cdot n)}$. Then $u (\kappa, \lambda) = \lambda$.
\item Suppose that $u (\kappa, \lambda) > \lambda$. Then $u (\kappa, \lambda) = u ((\kappa^{+(\kappa \cdot n)})^+, \lambda)$.
 \end{enumerate}
\end{fact}

\begin{Lem} Let $\kappa$ be a weakly inaccessible cardinal, and $\lambda$ be a cardinal such that $\kappa^{+(\kappa \cdot n)} < \lambda < \kappa^{+(\kappa \cdot (n+1))}$, where $n < \omega$, and $u (\kappa, \lambda) > \lambda$. Then ${\mathcal A}_{\kappa,\lambda} ((\cf (\theta))^+, u (\kappa, \lambda))$ holds, where $\theta = \theta ((\kappa^{+(\kappa \cdot n)})^+, \lambda)$. 
\end{Lem}

{\bf Proof.} By Observation 2.9 and Fact 3.8 (ii), it suffices to prove that ${\mathcal A}_{\kappa,\theta} ((\cf (\theta))^+, \pp (\theta))$ holds. By Fact 3.5, $\pp (\theta)$ is a successor cardinal, so we may find $A$ and $I$ such that
\begin{itemize}
\item $A$ is a set of regular cardinals smaller than $\theta$ with supremum $\theta$.
\item $\vert A \vert = \cf (\theta) < \min A$.
\item $I$ is an ideal on $A$ with $\{A \cap a : a \in A \} \subseteq I$.
\item  $\pp (\theta) = \tcf (\prod A /I )$.
\end{itemize}
Then by Fact 3.3 (iii), ${\mathcal A}_{\kappa,\theta} (\vert A \vert ^+, \pp (\theta))$ holds.
\hfill$\square$

\begin{Th} \begin{enumerate}[\rm (i)]
\item Let $\kappa$ be a weakly inaccessible cardinal, and $\lambda$ a cardinal with $\kappa \leq \lambda < \kappa^{+(\kappa \cdot \omega)}$. Then ${\mathcal A}_{\kappa,\lambda} (\tau, u (\kappa, \lambda))$ holds for some cardinal $\tau < \kappa$.
\item Let $\mu$ be a regular uncountable cardinal, $\kappa$ a weakly inaccessible cardinal greater than $\mu$, and $\lambda$ a cardinal with $\kappa \leq \lambda < \kappa^{+\mu}$. Then ${\mathcal A}_{\kappa,\lambda} (\mu, u (\kappa, \lambda))$ holds.
\end{enumerate}
\end{Th}

{\bf Proof.} By Fact 3.8 (ii) and Lemma 3.9.
\hfill$\square$

\bigskip

\subsection{C'me on, a little exponentiation will not harm you !}

\bigskip

The author would like to apologize to former addicts who would be trying hard to keep away from any kind of cardinal exponentiation. He respects their courageous fight and advises them to skip this section which might be painful to them.

\medskip

\begin{fact}  {\rm(\cite[Conclusion 1.6 (3) p. 321]{SheCA})}  Let $\nu$ and $\theta$ be two cardinals such that $\omega < \cf (\theta) \leq \nu < \theta$. Suppose that $\pp_\nu (\chi) < \theta$  for any large enough singular cardinal $\chi < \theta$ with $\cf(\chi) \leq \nu$. Then $\pp_\nu (\theta) = \pp(\theta)  = \pp_{\Gamma((\cf(\theta)^+, \cf(\theta))} (\theta)$.
\end{fact}

\begin{Obs} Let $\kappa$ be a regular cardinal, and $\theta$ be a singular cardinal such that $\tau^{< \kappa} < \theta$ for any cardinal $\tau < \theta$. Suppose that either $\theta$ is not a fixed point of the aleph function, or $\cf (\theta) > \omega$. Then the following hold :
\begin{enumerate}[\rm (i)]
\item $u (\kappa, \theta)$ equals $\theta$ if $\cf (\theta) \geq \kappa$, and $\pp (\theta)$ otherwise. 
\item  $\pp(\theta)  = \pp_{\Gamma((\cf(\theta)^+, \cf(\theta))} (\theta)$.
\end{enumerate}
\end{Obs}	

{\bf Proof.} (i) : If $\cf (\theta) \geq \kappa$, then $\theta \leq u (\kappa, \theta) \leq \theta^{< \kappa} = \theta$. Next assume that $\cf (\theta) < \kappa$. Then we have that $\pp_{< \kappa} (\theta) = \pp (\theta)$, by Fact 3.5 (ii) if $\theta$ is not a fixed point of the aleph function, and by Fact 3.11 otherwise. By results of Shelah (see e.g. Theorems 9.1.2 and 9.1.3 in \cite{HSW}), $\pp_\nu (\theta) = \theta^\nu$ for every cardinal $\nu$ with $\cf (\theta) \leq \nu < \kappa$. It follows that $\pp_{< \kappa} (\theta) = \theta^{< \kappa}$. Furthermore since $2^{< \kappa} < \theta \leq u (\kappa, \theta)$, $\theta^{< \kappa} = \max \{2^{< \kappa}, u (\kappa, \theta)\} = u (\kappa, \theta)$.

\medskip

(ii) : The assertion is trivial if $\cf (\theta) = \omega$, and follows from Fact 3.11 otherwise.
 \hfill$\square$

\begin{fact}  {\rm (\cite{Eis})}  Let $\sigma$, $\rho$ and $\theta$ be three cardinals such that $\omega < \sigma = \cf (\sigma) \leq \cf (\theta) < \rho = \cf (\rho) < \theta$. Suppose that $\pp_{\Gamma(\rho, \sigma)} (\theta)$ is regular, and moreover $u (\rho, \tau) < \pp_{\Gamma(\rho, \sigma)} (\theta)$ for any cardinal $\tau$ with $\rho \leq \tau < \theta$. Then $\pp_{\Gamma(\rho, \sigma)} (\theta) \in \PP_{\Gamma(\rho, \sigma)} (\theta)$.
\end{fact}

\begin{fact}  {\rm (\cite{Weaksat})}  Let $\tau$, $\kappa$, $\theta$ and $\rho$ be four uncountable cardinals such that $\tau \leq \kappa < \theta < \rho$. Suppose that ${\mathcal A}_{\kappa,\theta} (\tau, \sigma)$ holds for every regular infinite cardinal $\sigma \leq \rho$. Then ${\mathcal A}_{\kappa,\theta} (\tau, \rho)$ holds.
\end{fact}

\begin{Th} \begin{enumerate}[\rm (i)]
\item Let $\kappa$ be an infinite successor cardinal, say $\kappa = \nu^+$, and $\theta$ be a singular cardinal such that $\tau^\nu < \theta$ for any cardinal $\tau < \theta$. Suppose that either $\theta$ is not a fixed point of the aleph function, or $\cf (\theta) > \omega$. Suppose further that either $(\max \{\rho (\kappa), \cf (\theta)\})^{+3} < \kappa$, or $\cf (\theta) \in \rho (\kappa) \setminus \{ \cf (\rho (\kappa)) \}$. Then ${\mathcal A}_{\kappa,\theta} ((\cf (\theta)) ^+, u (\kappa, \theta))$ holds. 
\item Let $\kappa$ be a weakly inaccessible cardinal, and $\theta$ be a singular cardinal such that $\tau^{< \kappa} < \theta$ for any cardinal $\tau < \theta$. Suppose that either $\theta$ is not a fixed point of the aleph function, or $\cf (\theta) > \omega$. Then ${\mathcal A}_{\kappa,\theta} ((\cf (\theta)) ^+, u (\kappa, \theta))$ holds. 
\end{enumerate}
\end{Th}

{\bf Proof.} (i) : We can assume that $u (\kappa, \theta) > \theta$, since otherwise the conclusion is trivial. Then by Observation 3.12, $\cf (\theta) < \kappa$. By Fact 3.14 and, again, Observation 3.12, it suffices to prove that ${\mathcal A}_{\kappa,\theta} ((\cf (\theta)) ^+, \sigma)$ holds for any regular cardinal $\sigma$ with $\theta < \sigma \leq \pp_{\Gamma((\cf(\theta)^+, \cf(\theta))} (\theta)$. Given such a $\sigma$, we may find by Fact 3.13 $A$, $I$ and $\pi$ such that
\begin{itemize}
\item $A$ is a set of regular cardinals smaller than $\theta$ with supremum $\theta$.
\item $\vert A \vert = \cf (\theta) < \min A$.
\item $I$ is a $\cf (\theta)$-complete ideal on $A$ with $\{A \cap a : a \in A \} \subseteq I$.
\item  $\sigma \leq \pi = \tcf (\prod A /I )$.
\end{itemize}
Select an increasing, cofinal sequence $\vec{f} = \langle f_\alpha : \alpha < \pi \rangle$ in $(\prod A, <_I)$. Notice that we are  in one of the following three cases :
\begin{itemize}
\item $(\max \{ \rho (\kappa), \vert A \vert \})^{+3} < \kappa$.
\item $\vert A \vert < \cf (\rho (\kappa))$.
\item $\cf (\rho (\kappa)) < \vert A \vert < \rho (\kappa)$ and $I$ is $(\cf (\rho (\kappa))^+$-complete.
\end{itemize}
Hence by Fact 3.3 (i), there is a closed unbounded subset $C_{\vec{f}}$ of $\pi$, consisting of infinite limit ordinals, with the property that any $\delta$ in $C_{\vec{f}}$ with $\cf (\delta) = \kappa$ is a good point for $\vec f$. By Fact 3.3 (ii), it follows that ${\mathcal A}_{\kappa,\lambda} (\vert A \vert)^+, \pi)$ (and hence ${\mathcal A}_{\kappa,\lambda} (\vert A \vert)^+, \sigma)$) holds.

\medskip

(ii) : The proof is a (simpler) variant of that of (i), with Fact 3.3 (iii) substituting for Fact 3.3 ((i) and (ii)).
\hfill$\square$

\bigskip

\subsection{Isomorphisms}

\bigskip

We are now in a position to state results on the existence of isomorphisms between various $\mu$-club filters.

\medskip

\begin{Pro} Suppose that one of the following holds :
\begin{enumerate}[\rm (1)]
\item $\kappa$ is a weakly inaccessible cardinal, and $\lambda$ a cardinal with $\kappa \leq \lambda < \kappa^{+(\kappa \cdot \omega)}$. 
\item $\kappa$ is a weakly inaccessible cardinal, and $\lambda$ a singular cardinal such that $\tau^{< \kappa} < \lambda$ for any cardinal $\tau < \lambda$, and either $\lambda$ is not a fixed point of the aleph function, or $\cf (\lambda) > \omega$. 
\item $\kappa$ is an infinite successor cardinal, and $\lambda$ a singular cardinal such that either $\lambda$ is not a fixed point of the aleph function, or $\cf (\lambda) > \omega$. Furthermore $\tau^{< \kappa} < \lambda$ for any cardinal $\tau < \lambda$, and either $(\max (\rho (\kappa), \cf (\lambda)))^{+3} < \kappa$, or $\cf (\lambda) \in \rho (\kappa) \setminus \{ \cf (\rho (\kappa)) \}$. 
\end{enumerate}
Then for any sufficiently large regular cardinal $\mu < \kappa$, the ideals of the form $N\mu$-$S_{\kappa,\sigma}$, $\lambda \leq \sigma \leq u (\kappa, \lambda)$, are all isomorphic.\end{Pro}

{\bf Proof.} By Observation 3.2 and Theorems 3.10 (i) and 3.15.
\hfill$\square$

\begin{Pro}  Let $\lambda$ be a singular cardinal that is not a fixed point of the aleph function, and $\mu$ a regular cardinal with $\cf (\lambda) < \mu < \lambda$. Then for any sufficiently large regular cardinal $\kappa < \lambda$ such that either $(\rho (\kappa))^{+ 3} < \kappa$, or $\cf (\rho (\kappa)) \not= \cf (\lambda)$, the ideals of the form $N\mu$-$S_{\kappa,\sigma}$, $\lambda \leq \sigma \leq u (\kappa, \lambda)$, are all isomorphic.
\end{Pro}

{\bf Proof.} By Observation 3.2 and Theorem 3.6.
\hfill$\square$

\medskip

Concerning (2) and (3) of Proposition 3.16, note that if $\kappa$ and $\lambda$ are two uncountable cardinals such that $\kappa = \cf (\kappa) < \lambda$, and $\tau^{< \kappa} < \lambda$ for any cardinal $\tau < \lambda$, then for any cardinal $\sigma$ with $\lambda \leq \sigma \leq u (\kappa, \lambda)$, $u (\kappa, \sigma) = \sigma^{< \kappa} = \lambda^{< \kappa} = u (\kappa, \lambda)$ (and hence if $N\mu$-$S_{\kappa,\sigma}$ and $N\mu$-$S_{\kappa, u (\kappa, \lambda)}$ are isomorphic, then so are $N\mu$-$S_{\kappa,\sigma}$ and $N\mu$-$S_{\kappa, u (\kappa, \sigma)}$). Using the following, a similar remark can be made concerning Proposition 3.17.

\medskip

\begin{fact} {\rm (\cite{pcf})} Let $\kappa$ be an infinite successor cardinal, and $\lambda$ be a cardinal with $\kappa \leq \lambda < FP(\kappa)$. Then $u (\kappa, \sigma) = u (\kappa, \lambda)$ for any cardinal $\sigma$ with $\lambda \leq \sigma \leq u (\kappa, \lambda)$. 
\end{fact}

\begin{Pro} Suppose that one of the following holds :
\begin{enumerate}[\rm (1)] 
\item $\mu$ is a regular uncountable cardinal, $\kappa$ a successor cardinal, and $\lambda$ a cardinal with $\kappa < \lambda < \kappa^{+ \mu}$. Furthermore either $\mu \leq (\rho (\kappa))^+$ and $(\rho (\kappa))^{+ 3} < \kappa$, or $\mu \leq \cf (\rho (\kappa))$.
\item $\mu$ is a regular uncountable cardinal, $\kappa$ a weakly inaccessible cardinal greater than $\mu$, and $\lambda$ a cardinal with $\kappa \leq \lambda < \kappa^{+\mu}$. \end{enumerate}
Then the ideals of the form $N\mu$-$S_{\kappa,\sigma}$, $\lambda \leq \sigma \leq u (\kappa, \lambda)$, are all isomorphic. 
\end{Pro}

{\bf Proof.} By Observation 3.2 and Theorems 3.7 and 3.10 (ii).
\hfill$\square$

\bigskip

\section{Small $\mu$}

\bigskip

This case is more difficult. One reason for this should be clear. Suppose that $u (\kappa, \lambda) = \pp (\theta)$ for some $\theta$ about which we only know that $\cf (\theta) < \kappa < \theta \leq \lambda$. Then it cannot be ruled out that $\cf (\theta) = \omega$, in which case there is no infinite cardinal $\mu$ at all with $\mu < \cf (\theta)$. This is why in what follows we concentrate on situations when we have better information concerning $\cf (\theta)$ (e.g. when $\theta = \lambda$).

\medskip

Let $\sigma < \kappa \leq \lambda \leq \pi$ be four infinite cardinals, where $\kappa$ is regular, and $I$ be an ideal on $\sigma$. 
An ${\mathcal F}_{\kappa,\lambda}^I (\sigma^+, \pi)${\it - sequence}
is a $(\sigma^+,\lambda,\pi)$-sequence $\vec y = \langle y_\beta : \beta < \pi \rangle$ for which one may find a sequence $\vec f = \langle f_\beta : \lambda \leq \beta < \pi \rangle$ such that :
\begin{itemize}
\item For each $\beta \in \pi \setminus \lambda$, $f_\beta$ is a function from $\sigma$ to $ \lambda$.
\item  For any $e \in [\pi \setminus \lambda]^\kappa$, there are $b \in [e]^\kappa$ and $g : b \rightarrow I$ such that $f_\alpha (i) \not= f_\beta (i)$  whenever $\alpha < \beta$ are in $b$ and $i$ is in $\sigma \setminus (g(\alpha) \cup g(\beta))$.
\item For $\lambda \leq \beta < \pi$, $y_\beta = ran (f_\beta)$.
\end{itemize}

\medskip

\begin{fact} {\rm (\cite{Heaven})}
\begin{enumerate}[\rm (i)]
\item Any ${\mathcal F}_{\kappa,\lambda}^I   (\sigma^+, \pi)$- sequence is an ${\mathcal A}_{\kappa,\lambda} (\sigma^+,\pi)$-sequence.
\item Suppose that there exists an ${\mathcal F}_{\kappa,\lambda}^I   (\sigma^+, \pi)$- sequence, where $I$ is a $\mu^+$-complete ideal on $\sigma$. Then $ \pi \leq \cov (\lambda, \kappa, \kappa, \mu^+) $. 
\end{enumerate}
\end{fact}

\begin{Obs} Let $\pi$ be a singular cardinal greater than $\lambda$ with the property that for any cardinal $\tau$ with $\lambda \leq \tau < \pi$, there exists an ${\mathcal F}_{\kappa,\lambda}^I   (\sigma^+, \tau)$- sequence. Then there exists an ${\mathcal F}_{\kappa,\lambda}^I   (\sigma^+, \pi)$- sequence.
\end{Obs}

{\bf Proof.} Set $\rho = \cf (\pi)$, and select an increasing, continuous sequence $\langle \pi_\xi : \xi < \rho \rangle$ of cardinals greater than $\lambda$ with supremum $\pi$. For $\xi < \rho$, put $c_\xi = \{\beta : \pi_\xi \leq \beta < \pi_{\xi + 1}\}$, and pick $f^\xi_\beta : \sigma \rightarrow \lambda$ for $\beta \in c_\xi$ with the property that for any $w \in [c_\xi]^\kappa$, there are $b \in [w]^\kappa$ and $g : b \rightarrow I$ such that $f^\xi_\alpha (i) \not= f^\xi_\beta (i)$  whenever $\alpha < \beta$ are in $b$ and $i$ is in $\sigma \setminus (g(\alpha) \cup g(\beta))$. We may find $h_\xi : \sigma \rightarrow \lambda$ for $\xi < \rho$ such that that for any $v \in [\rho]^\kappa$, there are $b \in [v]^\kappa$ and $g : b \rightarrow I$ such that $h_\alpha (i) \not= h_\beta (i)$  whenever $\alpha < \beta$ are in $b$ and $i$ is in $\sigma \setminus (g(\alpha) \cup g(\beta))$. Fix a bijection $j : \lambda \times \lambda \rightarrow \lambda$. For $\xi < \rho$ and $\beta \in c_\xi$, define $t_\beta : \sigma \rightarrow \lambda$ by $t_\beta (i) = j (f_\beta^\xi (i), h_\xi (i))$.

\medskip

Let now $e \in [\pi]^\kappa$ be given.

\medskip

Case when there is $\zeta < \pi$ such that $\vert e \cap c_\zeta \vert = \kappa$. There must be $b \in [e \cap c_\zeta]^\kappa$ and $g : b \rightarrow I$ such that $f^\zeta_\alpha (i) \not= f^\zeta_\beta (i)$  whenever $\alpha < \beta$ are in $b$ and $i$ is in $\sigma \setminus (g(\alpha) \cup g(\beta))$. Then clearly, $t_\alpha (i) \not= t_\beta (i)$  whenever $\alpha < \beta$ are in $b$ and $i$ is in $\sigma \setminus (g(\alpha) \cup g(\beta))$.

\medskip

Case when $\vert e \cap c_\xi \vert < \kappa$ for all $\xi < \rho$. Pick $v \in [\rho]^\kappa$ so that $\vert e \cap c_\xi \vert = 1$ for all $\xi \in v$. We may find $b \in [v]^\kappa$ and $g : b \rightarrow I$ such that $h_\alpha (i) \not= h_\beta (i)$  whenever $\alpha < \beta$ are in $b$ and $i$ is in $\sigma \setminus (g(\alpha) \cup g(\beta))$. It is simple to see that $t_\alpha (i) \not= t_\beta (i)$  whenever $\alpha < \beta$ are in $b$ and $i$ is in $\sigma \setminus (g(\alpha) \cup g(\beta))$.
\hfill$\square$

\begin{fact} {\rm (\cite{Heaven})} Let $\mu$ be a regular infinite cardinal less than $\kappa$, and $\vec y = \langle y_\beta : \beta < \pi\rangle$ be
an ${\mathcal F}_{\kappa,\lambda}^I   (\sigma^+, \pi)$- sequence, where $I$ is a $\mu^+$-complete ideal on $\sigma$, as witnessed by $\vec f = \langle f_\beta : \lambda \leq \beta < \pi \rangle$. Define $\chi : P_\kappa (\lambda) \rightarrow P (\pi)$ by 

\centerline{$\chi (a) = a \cup \{ \beta \in \pi \setminus \lambda : \{i < \sigma : f_\beta (i) \in a \} \in I^+\}$.}

Then the following hold :
\begin{enumerate}[\rm (i)]
\item For any $b \in P_\kappa (\lambda)$, $\vert \chi (b) \vert < \kappa$ and $f_{\vec y} (b) \subseteq \chi (b)$.
\item For any $u : \mu \rightarrow P_\kappa (\lambda)$,  $\chi (\bigcup_{r < \mu} u (r)) \subseteq \bigcup_{r < \mu} \chi (u (r))$.
\item $\chi$ is a $(\mu, \kappa, \lambda, \pi)$-shuttle.
\end{enumerate}
\end{fact}

Thus if there exists an ${\mathcal F}_{\kappa,\lambda}^I   (\sigma^+, \pi)$- sequence, where $I$ is $\mu^+$-complete, then by Fact 2.6,  $N\mu$-$S_{\kappa,\lambda}$ and $N\mu$-$S_{\kappa, \pi}$ are isomorphic.

\medskip

\begin{fact} {\rm (\cite{Heaven})} Suppose that $\tcf (\prod A/I) = \pi$, where $A$ is an infinite set of regular infinite cardinals with $\vert A \vert = \sigma < \min (A)$ and $\kappa < \sup A \leq \lambda < \pi$, and $I$ an ideal on $A$, and that $\vec{f} = \langle f_\xi : \xi < \pi \rangle$ is an $<_I$-increasing, cofinal sequence in $(\prod A, <_I)$. Suppose further that for some closed unbounded subset $C$ of $\pi$, every $\delta$ in $C \cap E^\pi_\kappa$  is a good point for $\vec f$. Then there exists an ${\mathcal F}_{\kappa,\lambda}^I (\sigma^+, \pi)$- sequence.
\end{fact}

\begin{Th} Let $\theta$ be a singular cardinal that is not a fixed point of the aleph function. Then for any sufficiently large regular cardinal $\kappa < \theta$ such that either $(\rho (\kappa))^{+ 3} < \kappa$, or $\cf (\rho (\kappa)) \not= \cf (\theta)$, there exists an ${\mathcal F}_{\kappa,\lambda}^I   ((\cf (\theta))^+, u(\kappa, \theta))$- sequence, where $I$ is $\cf (\theta)$-complete. 
\end{Th}

{\bf Proof.} By Fact 3.5, $\pp (\theta)$ is a successor cardinal, so we may find $A$ and $I$ such that
\begin{itemize}
\item $A$ is a set of regular cardinals smaller than $\theta$ but greater than $\cf (\theta)$.
\item $\sup A = \theta$ and $\vert A \vert = \cf (\theta) < \min A$.
\item $I$ is a $\cf (\theta)$-complete ideal on $A$ with $\{A \cap a : a \in A \} \subseteq I$.
\item  $\pp (\theta) = \tcf (\prod A /I )$.
\end{itemize}
Select an increasing, cofinal sequence $\vec{f} = \langle f_\alpha : \alpha < \pp (\theta) \rangle$ in $(\prod A, <_I)$.

\medskip

By Fact 3.5, there must be $\eta$ with $\cf (\theta) < \eta < \theta$ such that for any regular cardinal $\kappa$ with $\eta \leq \kappa < \theta$, the following hold :
\begin{itemize}
\item $\kappa$ is not a fixed point of the aleph function.
\item $u (\kappa, \theta) > \theta = \theta (\kappa, \theta)$.
\end{itemize}
Fix a regular cardinal $\kappa$ such that $\eta \leq \kappa < \theta$, $\cf (\theta) \leq \rho (\kappa)$ and either $(\rho (\kappa))^{+ 3} < \kappa$, or $\cf (\rho (\kappa)) \not= \cf (\theta)$. By Fact 3.3 ((i) and (ii)), there is a closed unbounded subset $C$ of $\pp (\theta)$ such that every $\delta \in C$ of cofinality $\kappa$ is a good point for ${\vec f}$, and hence by Fact 4.4, there must exist an ${\mathcal F}_{\kappa,\lambda}^I   ((\cf (\theta))^+, \pp (\theta))$- sequence.
\hfill$\square$

\begin{Th} \begin{enumerate}[\rm (i)]
\item Let $\kappa$ be an infinite successor cardinal, and $\theta$ be a singular cardinal such that $\tau^{< \kappa} < \theta$ for any cardinal $\tau < \theta$. Suppose that either $\theta$ is not a fixed point of the aleph function, or $\cf (\theta) > \omega$ and $\pp (\theta)$ is regular. Suppose further that either $(\max \{\rho (\kappa), \cf (\theta)\})^{+3} < \kappa$, or $\cf (\theta) \in \rho (\kappa) \setminus \{ \cf (\rho (\kappa)) \}$. Then there exists an ${\mathcal F}_{\kappa,\lambda}^I   ((\cf (\theta))^+, u(\kappa, \theta))$- sequence, where $I$ is $\cf (\theta)$-complete. 
\item Let $\kappa$ be a weakly inaccessible cardinal, and $\theta$ be a singular cardinal such that $\tau^{< \kappa} < \theta$ for any cardinal $\tau < \theta$. Suppose that either $\theta$ is not a fixed point of the aleph function, or $\cf (\theta) > \omega$ and $\pp (\theta)$ is regular. Then there exists an ${\mathcal F}_{\kappa,\lambda}^I   ((\cf (\theta))^+, u(\kappa, \theta))$- sequence, where $I$ is $\cf (\theta)$-complete. 
\end{enumerate}
\end{Th}

{\bf Proof.} (i) : We can assume that $u (\kappa, \theta) > \theta$, since otherwise the conclusion is trivial. Then by Observation 3.11, $\cf (\theta) < \kappa$, and moreover $u (\kappa, \lambda) = \pp (\theta)$. Appealing again to Observation 3.11, and to Facts 3.5 and 3.12, we may find $A$ and $I$ such that
\begin{itemize}
\item $A$ is a set of regular cardinals smaller than $\theta$ but greater than $\cf (\theta)$.
\item $\sup A = \theta$ and $\vert A \vert = \cf (\theta)$.
\item $I$ is a $\cf (\theta)$-complete ideal on $A$ with $\{A \cap a : a \in A \} \subseteq I$.
\item  $\tcf (\prod A /I ) = \pp (\theta)$.
\end{itemize}
Select an increasing, cofinal sequence $\vec{f} = \langle f_\alpha : \alpha < \pi \rangle$ in $(\prod A, <_I)$. We must be in one of the following three cases :
\begin{itemize}
\item $(\max \{ \rho (\kappa), \vert A \vert \})^{+3} < \kappa$.
\item $\vert A \vert < \cf (\rho (\kappa))$.
\item $\cf (\rho (\kappa)) < \vert A \vert < \rho (\kappa)$ and $I$ is $(\cf (\rho (\kappa))^+$-complete.
\end{itemize}
Hence by Fact 3.3 (i), there is a closed unbounded subset $C_{\vec{f}}$ of $\pi$, consisting of infinite limit ordinals, with the property that any $\delta$ in $C_{\vec{f}}$ with $\cf (\delta) = \kappa$ is a good point for $\vec f$. By Fact 4.3, the existence of an ${\mathcal F}_{\kappa,\lambda}^I   ((\cf (\theta))^+, \pp (\theta))$- sequence follows.

\medskip

(ii) : The proof is a slight modification of that of (i).
\hfill$\square$

\medskip

Let us remark that (in spite of Observation 4.2) we do not know how to handle the case when $\theta$ is a fixed point of the aleph function of uncountable cofinality with $\pp (\theta)$ singular.

\medskip

We next list some consequences in terms of existence of isomorphisms. We will need the following that is readily checked.

\medskip

\begin{Obs} Let $\rho < \kappa \leq \lambda \leq \sigma \leq \pi$ be five infinite cardinals, where $\kappa$ is regular, and $I$ be an ideal on $\rho$. 
Suppose that there exists an ${\mathcal F}_{\kappa,\lambda}^I (\rho^+, \pi)$- sequence. Then there exists an ${\mathcal F}_{\kappa,\sigma}^I (\rho^+, \pi)$- sequence.
\end{Obs}

\begin{Pro} \begin{enumerate}[\rm (i)]
\item Let $\kappa$ be a weakly inaccessible cardinal, and $\lambda$ a singular cardinal of uncountable cofinality such that $\pp (\lambda)$ is regular. Suppose that $\tau^{< \kappa} < \lambda$ for any cardinal $\tau < \lambda$. Then for any regular cardinal $\mu < \cf (\lambda)$, the ideals of the form $N\mu$-$S_{\kappa,\sigma}$, $\lambda \leq \sigma \leq u (\kappa, \lambda)$, are all isomorphic.
\item Let $\kappa$ be an infinite successor cardinal, and $\lambda$ a singular cardinal of uncountable cofinality such that $\pp (\lambda)$ is regular. Suppose that $\tau^{< \kappa} < \lambda$ for any cardinal $\tau < \lambda$, and either $(\max (\rho (\kappa), \cf (\lambda)))^{+3} < \kappa$, or $\cf (\lambda) \in \rho (\kappa) \setminus \{ \cf (\rho (\kappa)) \}$. Then for any regular cardinal $\mu < \cf (\lambda)$, the ideals of the form $N\mu$-$S_{\kappa,\sigma}$, $\lambda \leq \sigma \leq u (\kappa, \lambda)$, are all isomorphic.
\item Let $\lambda$ be a singular cardinal of uncountable cofinality that is not a fixed point of the aleph function, and $\mu$ a regular cardinal less than $\cf (\lambda)$. Then for any sufficiently large regular cardinal $\kappa < \lambda$ such that either $(\rho (\kappa))^{+ 3} < \kappa$, or $\cf (\rho (\kappa)) \not= \cf (\lambda)$, the ideals of the form $N\mu$-$S_{\kappa,\sigma}$, $\lambda \leq \sigma \leq u (\kappa, \lambda)$, are all isomorphic.
\end{enumerate}
\end{Pro}

{\bf Proof.} (i) and (ii) : By Fact 4.2, Theorem 4.5 and Observation 4.6.

\medskip

(iii) : By Fact 4.2, Theorem 4.4 and Observation 4.6.

\hfill$\square$

\bigskip

\subsection{Nonsaturation}

\bigskip

How should we interpret the fact that $N\mu$-$S_{\kappa,\lambda}$ and $N\mu$-$S_{\kappa,u (\kappa, \lambda)}$ are isomorphic ? From a theoretical point of view, it can be seen as one more instance of a well-known phenomenon : if $\pi$ is the length of a scale on $\lambda$, that is $\pi = \tcf (\prod A/I)$ for some $A, I$ with $\sup A = \lambda$, then for any regular cardinal $\kappa$ with $\vert A \vert < \kappa < \lambda$, some properties of $P_\kappa (\pi)$ depend on those of $P_\kappa (\lambda)$. For applications the perspective is reversed, as we know much more about ideals on  $P_\kappa (\tau)$ in the case when $\tau$ is a regular cardinal. Thus when $N\mu$-$S_{\kappa,\lambda}$ and $N\mu$-$S_{\kappa,u (\kappa, \lambda)}$ are isomorphic (and $u (\kappa, \lambda)$ is regular), we see $N\mu$-$S_{\kappa,\lambda}$ as a copy of $N\mu$-$S_{\kappa,u (\kappa, \lambda)}$, a copy that inherits the properties of the original, in particular in terms of saturation. 

\medskip

We start with weak saturation, which is easier.

\medskip

Given an ideal $J$ on an infinite set $X$, a cardinal $\rho$ and $Y\subseteq P(X)$, $J$ is {\it $Y$-$\rho$-saturated} if there is no $Q \subseteq J^+$ with $\vert Q \vert = \rho$ such that $A \cap B \in Y$ for any two distinct members $A, B$ of $Q$.
 
We say that $J$ is {\it weakly $\rho$-saturated} if it is $\{\emptyset\}$-$\rho$-saturated, and {\it nowhere weakly $\rho$-saturated} if for any $B \in J^+$, $J \vert B$ is not weakly $\rho$-saturated. 

\medskip

\begin{fact} Let $\kappa$ be a regular uncountable cardinal, and $\pi$ a cardinal greater than or equal to $\kappa$. Then the following hold :
\begin{enumerate}[\rm (i)]
\item {\rm (\cite{Matsubara})} Suppose that $\kappa$ is a successor cardinal. Then no $\kappa$-complete, fine ideal on $P_\kappa (\pi)$ is weakly $\pi$-saturated.
\item {\rm (\cite{Weaksat})} Suppose that $\kappa$ is a weakly inaccessible cardinal. Then no normal, fine ideal on $J$ $P_\kappa (\pi)$ with $\{b \in P_\kappa (\pi) : \vert b \vert = \vert b \cap \kappa \vert \} \in J^\ast$ is weakly $\pi$-saturated.
\item {\rm (\cite{Ideals})} $\{b \in P_\kappa (\pi) : \vert b \vert = \vert b \cap \kappa \vert \} \in N\mu$-$S_{\kappa,\lambda}^\ast$.
\end{enumerate}
\end{fact}

\begin{Pro} Let $\kappa$ be a regular uncountable cardinal, and $\lambda$ a cardinal greater than or equal to $\kappa$. Suppose that $N\mu$-$S_{\kappa,\lambda}$ and $N\mu$-$S_{\kappa,u (\kappa, \lambda)}$ are isomorphic. Then $N\mu$-$S_{\kappa,\lambda}$ is nowhere weakly $u (\kappa, \lambda)$-saturated.
\end{Pro}

{\bf Proof.} By Fact 4.9, $N\mu$-$S_{\kappa,u (\kappa, \lambda)}$ is nowhere weakly $u (\kappa, \lambda)$-saturated. By Observation 2.2 (iv), the desired conclusion follows.
\hfill$\square$

\medskip

We now turn to ordinary saturation. Instead of large almost disjoint families, we will produce long towers.

\medskip

Given an ideal $J$ on an infinite set $X$, $Y\subseteq P (X)$ and a cardinal $\rho$, a $(J, Y)${\it -tower of length} $\rho$ is a sequence $\langle A_\alpha : \alpha < \rho \rangle$ of subsets of $X$ such that $(A_\alpha \setminus A_\beta, A_\beta \setminus A_\alpha) \in J^+ \times Y$ whenever $\alpha < \beta < \rho$.

\medskip

It is simple to see that if $\langle A_\alpha : \alpha < \rho \rangle$ is a $(J, Y)$-tower, then $\langle A_\alpha \setminus A_{\alpha + 1} : \alpha < \rho \rangle$ witnesses that $J$ is not $Y$-$\rho$-saturated.

\medskip

Let $\pi$ be a regular uncountable cardinal. The {\it bounding number} $\mathfrak{b}_\pi$ denotes the least cardinality of any $F \subseteq {}^\pi \pi$ with the property that there is no $g \in {}^\pi \pi$ such that $\vert \{\alpha < \pi : f (\alpha) \geq g (\alpha \} \vert < \pi$ for all $f \in F$. 

 



\medskip


\begin{fact}  {\rm (\cite{Heaven})} 
\begin{enumerate}[\rm (i)]
\item Let $\kappa, \lambda, \pi$ be three uncountable cardinals such that $\kappa = \cf (\kappa) \leq \lambda < \pi = \cf (\pi)$, and $k : P_\kappa (\lambda) \rightarrow P_\kappa (\pi)$. Further let $K$ and $G$ be two ideals on $P_\kappa( \lambda)$, and $J$ and $H$ two ideals on $P_\kappa (\pi)$. Suppose that $k (K) \subseteq J$, $H \subseteq \chi (G)$ and for some cardinal $\rho$, there is a $(J, H)$-tower $\langle A_\alpha : \alpha < \rho \rangle$. Then $\langle k^{- 1} (A_\alpha) : \alpha < \rho \rangle$  is a $(K, G)$-tower.
\item Let $\kappa$ be a regular uncountable cardinal, and $\pi$ a regular cardinal greater than $\kappa$. Then there exists an $(N\mu$-$S_{\kappa,\pi}, I_{\kappa, \pi})$-tower of length $\mathfrak{b}_\pi$.
\end{enumerate}
\end{fact} 

\begin{Pro} Let $\kappa$ be a regular uncountable cardinal, and $\lambda$ a cardinal greater than $\kappa$. Suppose that $u (\kappa, \lambda)$ is a regular cardinal, and $N\mu$-$S_{\kappa,\lambda}$ and $N\mu$-$S_{\kappa,u (\kappa, \lambda)}$ are isomorphic. Then there exists an $(N\mu$-$S_{\kappa,\lambda}, N\mu$-$S_{\kappa,\lambda})$-tower of length $\mathfrak{b}_{u (\kappa, \lambda)}$.
\end{Pro}

{\bf Proof.} Case when $u (\kappa, \lambda) = \lambda$. : By Fact 4.11 (ii), there must be an $(N\mu$-$S_{\kappa,\lambda}, I_{\kappa,\lambda})$-tower of length $\mathfrak{b}_{u (\kappa, \lambda)}$.

\medskip

Case when $u (\kappa, \lambda) > \lambda$ : By Fact 4.11 (ii), we may find an $(N\mu$-$S_{\kappa,u (\kappa, \lambda)}, N\mu$-$S_{\kappa,u (\kappa, \lambda)})$-tower $\langle A_\alpha : \alpha < \mathfrak{b}_{u (\kappa, \lambda)} \rangle$. Then by Observation 2.2 (iv) and Fact 4.11 (i),$\langle \chi^{- 1} (A_\alpha) : \alpha < \delta \rangle$ is an $(N\mu$-$S_{\kappa,\lambda}, N\mu$-$S_{\kappa,\lambda})$-tower.
\hfill$\square$

\medskip

We can get a stronger conclusion from a stronger assumption.

\medskip

\begin{fact}  {\rm (\cite{Heaven})} Let $\kappa$ be a regular uncountable cardinal, $\lambda$ a cardinal greater than or equal to $\kappa$, and $\pi$ a regular cardinal greater than or equal to $\lambda$. Suppose that there exists a $(\mu, \kappa, \lambda, \pi)$-shuttle $\chi$. Then there exists an $(N\mu$-$S_{\kappa,\lambda}, I_{\kappa, \lambda})$-tower of length $\mathfrak{b}_\pi$.
\end{fact}

\bigskip

\section{For the happy few : more about paradise in heaven}

\bigskip

We see $L$ as the core of heaven. If $V = L$, and $\mu$, $\kappa$ and $\lambda$ are three infinite cardinals with $\mu = \cf (\mu) < \kappa = \cf (\kappa) \leq \lambda$, then, as shown in \cite{Heaven}, the $\mu$-club filters on $P_\kappa (\lambda)$ and $P_\kappa (\lambda^{< \kappa})$ are isomorphic if and only if $\cf (\lambda) \not= \mu$. 

\medskip

The further away you go from the center, the less structure (and so the more freedom) there is. It does not take long until you lose GCH, provisionally replaced with Shelah's Strong Hypothesis (SSH) that asserts that $\pp (\theta) = \theta^+$ for every singular cardinal $\theta$. The exact consistency strength of the failure of SSH is not known, but believed to be about the same as that of the failure of SCH (the Singular Cardinal Hypothesis) that has been shown by Gitik \cite{Gitik}) to be equiconsistent with the existence of a measurable cardinal $\chi$ of Mitchell order $\chi^{++}$. Weak forms of square go a longer way (the failure of $\square_\tau^\ast$ for $\tau$ singular is known to entail a Woodin cardinal (see \cite[p. 702]{Jech2})). 


\medskip

Assuming that SSH holds, and moreover $\square_\sigma^\ast$ holds for every singular cardinal $\sigma$, it is shown in \cite{Heaven} that  $N\mu$-$S_{\kappa,\lambda}$ and $N\mu$-$S_{\kappa, u (\kappa, \lambda)}$ are isomorphic whenever $\mu < \kappa < \lambda$ are three infinite cardinals such that $\mu = \cf (\mu) < \kappa = \cf (\kappa) < \lambda$, $u (\kappa, \lambda)$ is regular, and $\cf (\lambda) \not= \mu$.

\medskip

After a while you may start wondering whether you are still in heaven. How to find out ? If for instance you live under ZFC + $(\star \star)$, where $(\star \star)$ is a statement that is true in $L$ but not provable in ZFC, does it mean that you are in heaven ? Hell no, since otherwise the \say{paradise} described in Subsection 3.5 would be part of heaven, and not of earth. Now it is true that under GCH, \emph{every} singular cardinal is a strong limit cardinal, but it is also true that in ZFC, \emph{some} singular cardinals are strong limit cardinals, and Theorem 3.15 (i) will apply to these.

\medskip

Some of the conditions in this Theorem 3.15 (i) can be removed at no great cost. For the argument to go through, we need every scale to be good (meaning that its points are almost all (in the sense of the club filter) good). By Fact 3.3 (i), we already know that some points are good. What about the others ? For one thing, they are good as well in case the length of the scale is the successor of a regular cardinal. 

\medskip

\begin{fact} {\rm (\cite{Good})} Let $A$ be an infinite set of regular cardinals with $\vert A \vert < \min A$, $I$ be an ideal on $A$ with $\{A \cap a : a \in A\} \subseteq I$, and $\vec{f} = \langle f_\alpha : \alpha < \pi \rangle$ be an increasing, cofinal sequence in $(\prod A, <_I)$, where $\pi$ is a cardinal greater than $\sup A$. Suppose that $\pi$ is the successor of a regular cardinal. Then there is a closed unbounded subset $C$ of $\pi$ with the property that any limit ordinal $\delta \in C$ with $\omega_1 \leq \cf (\delta) < \sup A$ is a good point for $\vec{f}$.
\end{fact}

\medskip

Let us next consider the case when the length of the scale is not the successor of a regular cardinal. 

\medskip

For two infinite cardinals $\theta$, $\pi$ with $\theta \leq \pi = \cf (\pi)$, the \emph{weakly approachable ideal} $I [\pi;\theta]$ is the collection of all $T \subseteq \pi$ such that for some closed unbounded subset $C$ of $\pi$, and some sequence $\langle u_\xi : \xi < \pi \rangle$ of elements of $P_\theta (\pi)$,
\begin{itemize}
\item each $\alpha$ in $T \cap C$ is an infinite limit ordinal with $\cf (\alpha) < \alpha$ ;
\item for any $\alpha \in T \cap C$, there is a cofinal subset $B$ of $\alpha$ of order-type $\cf (\alpha)$ such that $\{B \cap \eta : \eta < \alpha\} \subseteq \bigcup_{\xi < \alpha} P (u_\xi)$.
\end{itemize}

\medskip

\begin{fact} {\rm (\cite{Scales})} 
\begin{enumerate}[\rm (i)]
\item $I [\pi;\sigma] \subseteq I [\pi;\theta]$ for every infinite cardinal $\sigma \leq \theta$.
\item Suppose that $ \pi \notin I [\pi;\theta]$. Then $I [\pi;\theta]$ is a normal ideal on $\pi$.
\end{enumerate}
\end{fact}

\begin{fact} {\rm (\cite{Good})} Let $A$ be an infinite set of regular cardinals with $\vert A \vert < \min A$, $I$ be an ideal on $A$ with $\{A \cap a : a \in A\} \subseteq I$, and $\vec{f} = \langle f_\alpha : \alpha < \pi \rangle$ be an increasing, cofinal sequence in $(\prod A, <_I)$, where $\pi$ is a regular cardinal greater than $\sup A$. Further let $\kappa$ be a regular uncountable cardinal less than $\sup A$, and $S$ be a stationary subset of $E^\pi_\kappa$. Suppose that $S \in I [\pi;\sup A]$. Then there is a closed unbounded subset $C$ of $\pi$ with the property that any limit ordinal $\delta \in C \cap S$ is a good point for $\vec{f}$.
\end{fact}

\medskip

It is shown in \cite{Scales} that if $\pi$ is the successor of a cardinal $\tau$ and ${\rm WWS}_\pi (E^\pi_\kappa)$ (a consequence of $\square_\tau^\ast$) holds, then $E^\pi_\kappa \in I [\pi;\sup A]$.

\medskip

\begin{Th} \begin{enumerate}[\rm (i)]
\item Let $\theta$ be a singular cardinal that is not a fixed point of the aleph function. Suppose that either $\pp (\theta)$ is the successor of a regular cardinal, or $E^{\pp (\theta)}_\kappa \in I [\pp (\theta);\theta]$ for any sufficiently large regular cardinal $\kappa < \theta$ such that $\kappa = (\rho (\kappa))^{+ i}$ for $i = 1, 2$ or $3$, and $\cf (\rho (\kappa)) = \cf (\theta)$. Then for any sufficiently large regular cardinal $\kappa < \theta$, $u(\kappa, \theta) = \pp (\theta)$, and moreover ${\mathcal A}_{\kappa,\theta} ((\cf (\theta)) ^+, \pp (\theta))$ holds. 
\item Let $\kappa$ be an infinite successor cardinal, say $\kappa = \nu^+$, and $\theta$ be a singular cardinal such that $\tau^\nu < \theta$ for any cardinal $\tau < \theta$. Suppose that either $\theta$ is not a fixed point of the aleph function, or $\cf (\theta) > \omega$. Suppose further that one of the following holds :
\begin{itemize}
\item $\pp (\theta)$ is the successor of a regular cardinal.
\item $\pp (\theta)$ is regular, but not the successor of a regular cardinal, and if either $\kappa = (\cf (\theta))^{+ i}$ for $i = 1, 2$ or $3$, or $\cf (\theta) = \cf (\rho (\kappa)) < \rho (\kappa)$, then $E^{\pp (\theta)}_\kappa \in I [\pp (\theta);\theta]$.
\item $\pp (\theta)$ is singular, and if either $\kappa = (\cf (\theta))^{+ i}$ for $i = 1, 2$ or $3$, or $\cf (\theta) = \cf (\rho (\kappa)) < \rho (\kappa)$, then for cofinally many regular cardinals $\pi < \pp (\theta)$, $E^\pi_\kappa \in I [\pi;\theta]$.
\end{itemize}
Then ${\mathcal A}_{\kappa,\theta} ((\cf (\theta)) ^+, u (\kappa, \theta))$ holds. 
\end{enumerate}
\end{Th}

{\bf Proof.} (i) : By Fact 3.5, $\pp (\theta)$ is a successor cardinal, so there must be $A$ and $I$ such that
\begin{itemize}
\item $A$ is a set of regular cardinals smaller than $\theta$ with supremum $\theta$.
\item $\vert A \vert = \cf (\theta) < \min A$.
\item $I$ is an ideal on $A$ with $\{A \cap a : a \in A \} \subseteq I$.
\item $\pp (\theta) = \tcf (\prod A /I )$.
\end{itemize}


\medskip

Select an increasing, cofinal sequence $\vec{f} = \langle f_\alpha : \alpha < \pp (\theta) \rangle$ in $(\prod A, <_I)$.

\medskip

By Fact 3.5, we may find $\eta$ with $\cf (\theta) < \eta < \theta$ such that for any regular cardinal $\kappa$ with $\eta \leq \kappa < \theta$, the following hold :
\begin{itemize}
\item $\kappa$ is not a fixed point of the aleph function (so it is a successor cardinal).
\item $u (\kappa, \theta) = \pp (\theta)$.
\end{itemize}

\medskip

First assume that $\pp (\theta)$ is the successor of a regular cardinal. Then given a regular cardinal $\kappa$ with $\eta \leq \kappa < \theta$, there is by Fact 5.1 a closed unbounded subset $C$ of $\pp (\theta)$ with the property that any $\delta \in C \cap E^{\pp (\theta)}_\kappa$ is a good point for $\vec{f}$. By Fact 3.3, it follows that ${\mathcal A}_{\kappa,\theta} ((\cf (\theta)) ^+, \pp (\theta))$ holds. 

\medskip

Now assume that there is $\zeta$ with $\cf (\theta) < \zeta < \theta$ such that $E^{\pp (\theta)}_\kappa \in I [\pp (\theta);\theta]$ for any regular cardinal $\kappa$ with $\zeta \leq \kappa < \theta$ such that $\kappa = (\rho (\kappa))^{+ i}$ for $i = 1, 2$ or $3$, and $\cf (\rho (\kappa)) = \cf (\theta)$. Let $\kappa$ be a regular cardinal with $\max \{\eta, \zeta\} \leq \kappa < \theta$. If $\kappa = (\rho (\kappa))^{+ i}$ for $i = 1, 2$ or $3$, and $\cf (\rho (\kappa)) = \cf (\theta)$, then by Fact 5.3, there is a closed unbounded subset $C$ of $\pp (\theta)$ with the property that any $\delta \in C \cap E^{\pp (\theta)}_\kappa$ is a good point for $\vec{f}$, so by Fact 3.3, ${\mathcal A}_{\kappa,\theta} ((\cf (\theta)) ^+, \pp (\theta))$ holds.

Fix a regular cardinal $\kappa$ such that $\eta \leq \kappa < \theta$, $\cf (\theta) \leq \rho (\kappa)$ and either $(\rho (\kappa))^{+ 3} < \kappa$, or $\cf (\rho (\kappa)) \not= \cf (\theta)$. By Fact 3.3 ((i) and (ii)), there is a closed unbounded subset $C$ of $\pp (\theta)$ such that every $\delta \in C$ of cofinality $\kappa$ is a good point for ${\vec f}$, and hence ${\mathcal A}_{\kappa,\theta} (\vert A \vert ^+, \pp (\theta))$ holds.

\medskip

(ii) : The proof is a modification of that of Theorem 3.15 (i).
\hfill$\square$


\begin{Cor} 
\begin{enumerate}[\rm (i)]
\item Let $\theta$ be a singular cardinal that is not a fixed point of the aleph function. Suppose that $\tau^+ \in I [\tau^+;\theta]$ for any singular cardinal $\tau$ with $\theta \leq \tau < \pp (\theta)$. Then for any sufficiently large regular cardinal $\kappa < \theta$, $u(\kappa, \theta) = \pp (\theta)$, and moreover ${\mathcal A}_{\kappa,\theta} ((\cf (\theta)) ^+, \pp (\theta))$ holds. 
\item Let $\kappa$ be an infinite successor cardinal, say $\kappa = \nu^+$, and $\theta$ be a singular cardinal such that $\tau^\nu < \theta$ for any cardinal $\tau < \theta$. Suppose that either $\theta$ is not a fixed point of the aleph function, or $\cf (\theta) > \omega$ and $\pp (\theta)$ is not weakly inaccessible. Suppose further that $\tau^+ \in I [\tau^+;\theta]$ for any singular cardinal $\tau$ with $\theta \leq \tau < \pp (\theta)$. Then ${\mathcal A}_{\kappa,\theta} ((\cf (\theta)) ^+, u (\kappa, \theta))$ holds. 
\end{enumerate}
\end{Cor}

{\bf Proof.}  By Theorems 3.6, 3.15 (i) and 5.4.
\hfill$\square$

\begin{Th}  Let $\kappa$ be an infinite successor cardinal, and $\lambda$ a cardinal such that $\kappa < \lambda < \kappa^{+ (\rho (\kappa))^+}$ and $u (\kappa, \lambda) > \lambda$. Suppose that one of the following holds :
\begin{itemize}
\item $\pp (\theta (\kappa, \lambda))$ is the successor of a regular cardinal.
\item If $\kappa = (\rho (\kappa))^{+ i}$ for $i = 1, 2$ or $3$, then $E^{\pp (\theta (\kappa, \lambda))}_\kappa \in I [\pp (\theta (\kappa, \lambda));\theta (\kappa, \lambda)]$.
\end{itemize}
Then ${\mathcal A}_{\kappa,\lambda} ((\cf (\theta (\kappa, \lambda)) ^+, u (\kappa, \lambda))$ holds. 
\end{Th}

{\bf Proof.} The proof is similar to that of Theorem 5.4 (i).
\hfill$\square$


\begin{Cor} Let $\kappa$ be an infinite successor cardinal, and $\lambda$ a cardinal such that $\kappa < \lambda < \kappa^{+ (\rho (\kappa))^+}$. Suppose that $\tau^+ \in I [\tau^+;\theta]$ for any two singular cardinals $\tau$ and $\theta$ with $\kappa < \theta \leq \lambda \leq \tau < u (\kappa, \lambda)$. Then ${\mathcal A}_{\kappa,\lambda} ((\rho (\kappa))^+, u (\kappa, \lambda))$ holds. 
\end{Cor}

{\bf Proof.} Use Theorem 3.7 in case $(\rho (\kappa))^{+ 3} < \kappa$, and Theorem 5.6 otherwise.
\hfill$\square$

\bigskip

\section{Your choice of paradise}

\bigskip

Aficionados will not balance. They will tell you that yes, the heavenly paradise is (way) bigger, but its earthly cousin is free for all.
Besides, earth is more contrasted, and as theorized by the leibnizian doctor in Voltaire's Candide, \say{les malheurs particuliers font le bien g\'en\'eral ; de sorte que plus il y a de malheurs particuliers, et plus tout est bien}. In other words, how can one possibly enjoy oneself when everything is uniformly nice ? As the Eva of Klopstock's Messias puts it, \say{Z\"{a}rtlich seh'ich, mit irrendem Blick hinunter zur Erde ; Dich, Paradies, dich seh'ich nicht mehr}.

  \bigskip
\noindent Universit\'e de Caen - CNRS \\
Laboratoire de Math\'ematiques \\
BP 5186 \\
14032 Caen Cedex\\
France\\
Email :  pierre.matet@unicaen.fr\\


\end{document}